\documentclass[12pt]{article}
\usepackage{amssymb,amsmath,theorem}

\let\eps=\varepsilon \let\kappa=\varkappa 
\let\cal=\mathcal  \let\goth=\mathfrak  
\newcommand\norm[1]{\left\|#1\right\|}           
\newcommand\qed{\ifhmode\unskip\nobreak\fi\quad  
   \ifmmode\square\else\hbox{$\square$}\fi}      
\newcommand\proofskip{\vspace{
   \theorempostskipamount}}                      

\newtheorem{theorem}{Theorem}
\newtheorem{lemma}[theorem]{Lemma}

{\theorembodyfont{\normalfont}\newtheorem{rem}[theorem]{Remark}
                        
                        \newtheorem{exm}[theorem]{Example}
                        }

\newcommand\proof[1]{\noindent\textit{Proof#1}}

\begin{document}

\null\vspace{9pt}

\begin{center}
\Large\bfseries A new definition of $\boldsymbol{t}$-entropy\\ for transfer operators

\bigskip\bigskip\normalsize\rm

V.\,I.\ BAKHTIN

\smallskip

{\it Belarus State University, Belarus (e-mail: bakhtin@tut.by)}

\bigskip

A.\,V.\ LEBEDEV\footnote{Supported by the Grant of Polish Minister of Science and Higher Education
N201\,382634}

\smallskip

{\it University of Bialystok, Poland\ \ {\rm\&}\ \ Belarus State University, Belarus

(e-mail: lebedev@bsu.by)}
\end{center}

\renewcommand\abstractname{}
\begin{abstract}
The article presents a new definition of $t$-entropy that makes it more explicit and simplifies the
process of its calculation.
\end{abstract}

\bigbreak\bigskip

\quad\parbox{13.9cm} {\textbf{Keywords:} {\itshape transfer operator, variational principle,
$t$-entropy}

\medbreak \textbf{2000 MSC:} 37A35; 47B37; 47C15}

\vspace{10mm}

In the series of articles \cite{ABL-2000-1,ABL-2000-2,ABL-2001,ABL-2005,Bakh,ABL-2008} there have
been established the variational principles for the spectral radii  of weighted shift and transfer
operators generated by an arbitrary dynamical system. These principles are based on the Legendre
duality and the main role here is played by a newly introduced dynamical invariant --- $t$-entropy,
which gives the explicit form of the Legendre dual object to the logarithm of the spectral radii of
operators in question.  The description of \mbox{$t$-entropy} is not elementary and its calculation
is rather sophisticated. In the present article we give a new definition of $t$-entropy that makes
it more explicit and essentially simplifies the process of its calculation.

The article consists of two sections. In Section \ref{1..} we consider $t$-entropy for the   model
example of transfer operators associated with  continuous dynamical systems. The new definition of
$t$-entropy is introduced here in Theorem \ref{1..2}. In Section \ref{2} we discuss the general
$C^*$-dynamical situation. To illustrate similarity and difference between the objects considered
in the model and general situations we  present here  a number of  examples and finally introduce
the general new definition of $t$-entropy in Theorem \ref{2..7}.

\section{A new definition of $\boldsymbol{t}$-entropy for continuous\\ dynamical systems}\label{1..}

Let us consider a Hausdorff compact space $X$. We denote by  $C(X)$ the algebra of continuous
real-valued functions on~$X$ equipped with the uniform norm. Let  $\alpha\!:X\to X$ be a continuous
mapping. This mapping generates the dynamical system with  discrete time, which will be denoted by
$(X,\alpha)$.

A linear operator  $A\!:C(X)\to C(X)$ is called a \emph{transfer operator} for the dynamical system
$(X,\alpha)$ if

a) $A$  is positive  (that is it maps nonnegative functions to nonnegative) and

b) it satisfies the  \emph{homological identity}
\begin{equation}\label{1,,1}
A\bigl(f\circ \alpha\cdot g\bigr) = fAg, \qquad f,g\in C(X).
\end{equation}

We denote by  $M\subset C^*(X)$ the set of all linear positive normalized functionals on~$C(X)$
(that is linear functionals that take nonnegative values on nonnegative functions and are equal to
$1$ on the unit function). By the Riesz theorem these functionals are bijectively identified with
regular Borel probability measures on~$X$, and so all the elements of $M$ will be referred to as
measures.

A measure  $\mu\in M$ is called  $\alpha$-invariant if  $\mu(f\circ\alpha)  = \mu(f)$ for all
functions $f\in C(X)$. This is equivalent to the identity $\mu(\alpha^{-1}(G)) = \mu(G)$ for all
Borel subsets  $G\subset X$. The collection of all $\alpha$-invariant measures from $M$ will be
denoted by $M_\alpha$.

Recall that a continuous partition of unity on $X$ is a finite set $D = \{g_1,\dotsс,g_k\}$
consisting of nonnegative functions $g_i\in C(X)$ satisfying the identity $g_1+ \dots +g_k \equiv
1$.

According to \cite{ABL-2008}, \emph{$t$-entropy} is the  functional $\tau(\mu)$ on $M$ defined by
the formulae
\begin{gather}\label{1,,2}
 \tau(\mu) :=\inf_{n\in\mathbb N}\frac{\tau_n(\mu)}{n}\,,\qquad
 \tau_n(\mu) :=\inf_D\tau_n(\mu,D),\\[6pt]
 \tau_n(\mu,D) :=\sup_{m\in M}\sum_{g\in D}\mu(g) \ln\frac{m(A^ng)}{\mu(g)}\,. \label{1,,3}
\end{gather}
The second infimum in \eqref{1,,2} is taken over all the continuous partitions of unity $D$ on~$X$.
If we have $\mu(g)  = 0$ for a certain element $g\in D$, then we set the corresponding summand
in~\eqref{1,,3} to be zero independently of the value~$m(A^ng)$. And if there exists an element
$g\in D$ such that $A^ng =  0$ and simultaneously $\mu(g)>0$, then we set $\tau(\mu)  = -\infty$.

Given a transfer operator $A$ we define a family of operators $A_\varphi\!:C(X) \to C(X)$ depending
on the functional parameter $\varphi\in C(X)$ by means of the formula
$$
 A_\varphi f =A(e^\varphi f).
$$
Evidently, all the operators of this family are transfer operators as well. Let us denote
by~$\lambda(\varphi)$ the logarithm of the spectral radius of $A_\varphi$, that is
$$
\lambda(\varphi)  = \lim_{n\to\infty}\frac{1}{n}\ln \norm{A_\varphi^n}.
$$

The principal importance of $t$-entropy is clearly demonstrated by the following Variational
Principle.

\begin{theorem}\label{1..1}{\rm (\cite{ABL-2008}, Theorem 5.6)\,}
Let\/ $A\!:C(X)\to C(X)$ be a transfer operator for a continuous mapping\/  $\alpha\!:X\to X$ of a
Hausdorff compact space\/~$X$.  Then
\begin{equation*}
\lambda(\varphi)= \max_{\mu\in M_\alpha} \bigl(\mu(\varphi)+\tau(\mu)\bigr),\qquad\varphi\in C(X).
\end{equation*}
\end{theorem}

The next theorem presents a new definition of  $t$-entropy.

\begin{theorem}\label{1..2}
For $\alpha$-invariant measures\/ $\mu\in M_\alpha$ the following formula is true
\begin{equation}\label{1,,4}
\tau(\mu)  =\inf_{n,D} \frac{1}{n} \sum_{g\in D}\mu(g) \ln\frac{\mu(A^ng)}{\mu(g)}.
\end{equation}
\end{theorem}

In other words, in the definition of $t$-entropy for an $\alpha$-invariant measure $\mu$ one should
not calculate the supremum in \eqref{1,,3} but can simply put $m=\mu$ there. Thus, expression
\eqref{1,,3} is changed for
\begin{equation}\label{1,,5}
\tau'_n(\mu,D) =\sum_{g\in D}\mu(g) \ln\frac{\mu(A^ng)}{\mu(g)}.
\end{equation}

\medskip

To prove Theorem \ref{1..2} we need the next

\begin{lemma}\label{1..3}
For any continuous partition of unity\/ $D$ on\/ $X$ and any pair of numbers\/ $n\in \mathbb N$,\ \
$\eps>0$ there exists a continuous partition of unity\/ $E$ such that for each pair of functions\/
$g\in D$ and\/ $h\in E$ the oscillation of\/ $A^ng$ on the support of\/ $h$ is less than\/ $\eps$:
\begin{equation}\label{1,,6}
\sup \bigl\{A^ng(x)\bigm| h(x)>0\bigr\} -\inf\bigl\{A^ng(x)\bigm| h(x)>0\bigr\} <\eps.
\end{equation}
\end{lemma}

\smallskip

\proof. Any function $A^ng$ belongs to $C(X)$. Therefore its range is contained in a certain
segment $[a,b]$.

Evidently, there exists a continuous partition of unity  $\{f_1,\dots,f_k\}$ on the segment $[a,b]$
such that the support of every its element is contained in a certain interval of the length less
than $\eps$. Then the family $E_g =\{f_1\circ\! A^ng,\,\ldots,\,f_k\circ\! A^ng\}$ forms a
continuous partition of unity on $X$ and on the support of each its element the oscillation of
$A^ng$ is less than~$\eps$. Now all the products $\prod_{g\in D} h_g$, where $h_g\in E_g$, form the
desired partition of unity~$E$. \qed \proofskip

Now let us prove Theorem  \ref{1..2}. Comparing  \eqref{1,,3} and  \eqref{1,,5} one sees that
$$\tau'_n(\mu,D) \le \tau_n(\mu,D).$$
Therefore to prove  \eqref{1,,4} it is enough to verify the inequality
$$
\tau_n(\mu)\le\tau'_n(\mu,D).
$$
Since in the case when  $\tau_n(\mu) = -\infty$ the latter inequality is trivial in what follows we
assume that  $\tau_n(\mu) >-\infty$.

Let us fix a number $n\in\mathbb N$, a continuous partition of unity $D$ on $X$ and a number
$\eps>0$. For these objects there exists a continuous partition of unity $E$ mentioned in Lemma
\ref{1..3}. Consider one more continuous partition of unity consisting of the functions of the form
$g\cdot h\circ\alpha^n$, where $g\in D$ and  $h\in E$. For this partition, by the definition of
$\tau_n(\mu)$, there exists a measure  $m\in M$ such that
\begin{equation*}
\tau_n(\mu)-\eps\le\sum_{g\in D}\sum_{h\in E} \mu(g\cdot h\circ\alpha^n)\ln \frac{m\bigl(A^n(g\cdot
h\circ\alpha^n) \bigr)}{\mu(g\cdot h\circ\alpha^n)}.
\end{equation*}
From the homological identity it follows that $A^n(g\cdot h\circ\alpha^n) =hA^ng$. Therefore, the
latter inequality is equivalent to
\begin{equation}\label{1,,7}
\tau_n(\mu)-\eps\le\sum_{g\in D}\sum_{h\in E} \mu(g\cdot h\circ\alpha^n)\ln
\frac{m(hA^n(g))}{\mu(g\cdot h\circ\alpha^n)}.
\end{equation}

Now for each pair $g\in D$,\,\ $h\in E$ choose a number $y_{gh}$ satisfying two conditions
\begin{gather}\label{1,,8}
m(hA^ng) =m(h)y_{gh},\\[6pt] \label{1,,9}
\inf\bigl\{A^ng(x)\bigm|h(x)>0\bigr\} \le y_{gh}\le \sup\bigl\{A^ng(x)\bigm|h(x)>0\bigr\}.
\end{gather}
Then inequality \eqref{1,,7} takes the form
\begin{equation}\label{1,,10}
\tau_n(\mu)-\eps\le\sum_{g\in D}\sum_{h\in E} \mu(g\cdot h\circ\alpha^n)\ln
\frac{m(h)y_{gh}}{\mu(g\cdot h\circ\alpha^n)}\,,
\end{equation}
which is equivalent to
\begin{equation}\label{1,,11}
\tau_n(\mu)-\eps\le
 \sum_{g\in D}\sum_{h\in E} \mu(g\cdot h\circ\alpha^n)\ln \frac{y_{gh}}{\mu(g\cdot h\circ\alpha^n)}
 + \sum_{g\in D}\sum_{h\in E} \mu(g\cdot h\circ\alpha^n)\ln {m(h)}.
\end{equation}

Let us consider separately the second summand in the right-hand side of \eqref{1,,11}:
\begin{equation}\label{1,,12}
\sum_{g\in D}\sum_{h\in E} \mu(g\cdot h\circ\alpha^n)\ln m(h) =\sum_{h\in E} \mu(h\circ
\alpha^n)\ln m(h) =\sum_{h\in E} \mu(h)\ln m(h).
\end{equation}
Here in the left-hand equality we have exploited the fact that $D$ is a partition of unity and in
the right-hand equality we have used $\alpha$-invariance of $\mu$. If we treat $m(h)$ in
\eqref{1,,12} as independent nonnegative variables satisfying the condition $\sum_{h\in E}m(h) =1$
then the routine usage of Lagrange multipliers principle shows that the function $\sum_{h\in E}
\mu(h)\ln m(h)$ attains its maximum when $m(h) =\mu(h)$. Evidently, the same is true for the
right-hand sides in \eqref{1,,11} and \eqref{1,,10}. Therefore,
\begin{equation}\label{1,,13}
\tau_n(\mu)-\eps\le\sum_{g\in D}\sum_{h\in E} \mu(g\cdot h\circ\alpha^n)\ln
\frac{\mu(h)y_{gh}}{\mu(g\cdot h\circ\alpha^n)}\,.
\end{equation}

Observe that estimates \eqref{1,,6} and \eqref{1,,9} imply
\begin{equation}\label{1,,14}
\mu(h)y_{gh} \le \mu\bigl(h(A^ng +\eps)\bigr).
\end{equation}
Using \eqref{1,,13}, \eqref{1,,14}, and the fact that $E$ is a partition of unity and exploiting
the concavity of logarithm we obtain the following relations:
\begin{align*}
\tau_n(\mu)-\eps & \le  \sum_{g\in D}\sum_{h\in E} {\mu(g\cdot h\circ\alpha^n)}
\ln\frac{\mu\bigl(h(A^ng +\eps)\bigr)}{\mu(g\cdot h\circ\alpha^n)} =\\[3pt]
&=\sum_{g\in D}\mu(g)\sum_{h\in E} \frac{\mu(g\cdot h\circ\alpha^n)}{\mu(g)} \ln
\frac{\mu\bigl(h(A^ng +\eps)\bigr)}{\mu(g\cdot h\circ\alpha^n)} \le\\[3pt]
&\le\sum_{g\in D}\mu(g)\ln\sum_{h\in E} \frac{\mu\bigl(h(A^ng +\eps)\bigr)}{\mu(g)} = \sum_{g\in D}
\mu(g)\ln\frac{\mu(A^ng+\eps)}{\mu(g)}.
\end{align*}
By the arbitrariness of $\eps$ this implies
\begin{equation*}
 \tau_n(\mu)\le \sum_{g\in D}\mu(g)\ln\frac{\mu(A^ng)}{\mu(g)} =\tau'_n(\mu,D)
\end{equation*}
and finishes the proof of Theorem~\ref{1..2}. \qed

\proofskip

Now let us proceed to the general $C^*$-dynamical setting.

\section{The general case of $\boldsymbol{C^*}$-dynamical systems}\label{2}

The general notion of \mbox{$t$-entropy} involves the so-called base algebra and a transfer
operator for a $C^*$-dynamical system. Let us recall definitions of these objects (see
\cite{ABL-2008}).

A \emph{base algebra} $\cal C$ is a selfadjoint part of a certain commutative $C^*$-algebra with an
identity $\bf 1$. This means that there exists a commutative $C^*$-algebra $\cal B$ with an
identity $\bf 1$ such that
$$
{\cal C} = \{\,b\in {\cal B}\mid b^* = b\,\}.
$$

\medskip

A pair $({\cal C},\delta)$, where $\cal C$ is a base algebra and $\delta$ is its certain
endomorphism such that $\delta({\mathbf 1}) = {\mathbf 1}$, is called a $C^*$-\emph{dynamical
system}.

Let $({\cal C},\delta)$ be a $C^*$-dynamical system. A linear operator $A\!:\cal C\to \cal C$ is
called a \emph{transfer operator} (for $(\cal C,\delta)$), if it possesses the following two
properties

a) $A$ is positive (it maps nonnegative elements  of $\cal C$ into nonnegative ones);

b) it satisfies the  \emph{homological identity}
\begin{equation}\label{2,,1}
 A\bigl((\delta f)g\bigr) =fAg \quad \textrm{for all}\ \ f,g\in \cal C.
\end{equation}

We denote by $M(\cal C)$ the set of all positive normalized linear functionals on~$\cal C$. A
functional $\mu\in M(\cal C)$ is called  \emph{$\delta$-invariant} if for each $f\in \cal C$ we
have $\mu(f)=\mu(\delta f)$. The set of all $\delta$-in\-va\-ri\-ant functionals from $M(\cal C)$
will be denoted by $M_\delta(\cal C)$.

By a \emph{partition of unity} in the algebra $\cal C$ we mean any finite set $D =
\{g_1,\dotsс,g_k\}$ consisting of nonnegative elements $g_i\in \cal C$ satisfying the identity
$g_1+ \dots +g_k = {\mathbf 1}$.

The definition of $t$-entropy introduced in the previous section in \eqref{1,,2} and \eqref{1,,3}
can be carried over word by word to the case of $C^*$-dynamical systems. Namely, here
\emph{$t$-entropy} is the  functional $\tau$  on $M(\cal C)$ such that its value at any $\mu\in
M(\cal C)$ is defined by the following formulae
\begin{gather}\label{2,,2}
 \tau(\mu)  := \inf_{n\in\mathbb N}\frac{\tau_n(\mu)}{n}\,,\qquad
 \tau_n(\mu)  := \inf_D\tau_n(\mu,D),\\[6pt]
 \tau_n(\mu,D)  := \sup_{m\in M (\cal C)}\sum_{g\in D}\mu(g) \ln\frac{m(A^ng)}{\mu(g)}\,. \label{2,,3}
\end{gather}
The infimum in \eqref{2,,2} is taken over all the partitions of unity $D$ in the algebra $\cal C$.

This $t$-entropy plays a crucial role in the  corresponding variational principles for the spectral
radii as for abstract transfer operators, so also for weighted shift operators in \hbox{$L^p$-type}
spaces (see \cite{ABL-2008}, Theorems 6.10, 11.2, 13.1 and 13.6).

To illustrate similarity and difference between the objects considered in this and the previous
sections we present now a number of examples of $C^*$-dynamical systems and transfer operators that
show, in particular, how far away from the continuous setting described in Section \ref{1..} one
can move.

\begin{exm}\label{2..1}
Let $Y$ be a measurable space with a $\sigma$-algebra $\goth A$ and $\beta\!:Y\to Y$ be a
measurable mapping. We denote by $(Y,\beta)$ the discrete time dynamical system generated by the
mapping $\beta$ on the phase space $Y$. Let $\cal C$ be any Banach algebra such that

a) $\cal C$ consists of bounded real-valued measurable functions on $Y$,

b) it is supplied with the uniform norm,

c) it contains the unit function, and

d) it is $\beta$-invariant (that is $f\circ\beta\in\cal C$ for all  $f\in\cal C$).

\smallskip

Clearly the mapping $\delta\!: \cal C\to \cal C$ given by $\delta(f) := f\circ\beta$ is an
endomorphism of $\cal C$ and therefore $(\cal C, \delta)$ is a $C^*$-dynamical system with the base
algebra $\cal C$.
\end{exm}

\begin{exm}\label{2..2}
As a particular case of the base algebra in the previous example one can take the algebra  of all
bounded real-valued measurable functions on $Y$. We will denote this algebra by~$B(Y)$.
\end{exm}

\begin{exm}\label{2..3}
Let $(Y,\goth A,m)$ be a measurable space with a probability measure $m$, and let~$\beta$ be a
measurable mapping such that $m\bigl(\beta^{-1}(G)\bigr)\le Cm(G),\ \ G\in \goth A$, where the
constant $C$ does not depend on $G$. In this case one can take as a base algebra the space
$L^\infty(Y,m)$ of all essentially bounded real-valued measurable functions on $Y$ with the
essential supremum norm.
\end{exm}

\begin{rem}\label{2..4}
1) If, as  in Example \ref{2..2},\ \ $\cal C =B(Y)$ then the elements of~$M(\cal C)$ can be
naturally identified with finitely-additive probability measures on the $\sigma$-algebra~$\goth A$
by means of the equality $\mu(f) =\int_Y f\,d\mu,\ \ f\in B(Y)$.

2) If, as in Example \ref{2..3},\ \ $\cal C = L^\infty (Y,m)$ then $M(\cal C)$ consists of
finitely-additive probability measures on $\goth A$ which are absolutely continuous with respect
to~$m$ (that is they are equal to zero on the sets of zero measure $m$).

3) In  Example \ref{2..3} the set $M_\delta (\cal C)$ is the subset of $M(\cal C)$ consisting of
measures $\mu$ such that $\mu(\beta^{-1}(G)) =\mu(G)$ for each measurable set~$G$.

4) It should be emphasized  that in general given a specific functional algebra its endomorphism is
\emph{not} necessarily generated by a \emph{point mapping} of the domain. For example, if $\cal C =
L^\infty(Y,m)$ then its endomorphisms are generated by \emph{set mappings} that do not `feel' sets
of measure zero
(see, for example, \cite{Walt2}, Chapter~2). Thus not every endomorphism of \ $ L^\infty (Y,m)$ is
generated by a certain measurable mapping~$\beta$ as in Example \ref{2..3}.

On the other hand, on the maximal ideals level any endomorphism is induced by a certain point
mapping (for details see \cite{ABL-2008}, Theorem~6.2). Therefore raising the apparatus of
investigation  to the $C^*$-algebraic level we not only essentially extend the field of its
applications but additionally can always exploit point mappings in the study of transfer operators
independently of their concrete origin (see in this connection the general description of transfer
operators given in \cite{ABL-2008}, Section~7).
\end{rem}

The next example can be considered as a model example of transfer operators on $L^\infty$.

\begin{exm}\label{2..5}
Let $(Y,\goth A)$ be a measurable space with a $\sigma$-finite measure $m$, and let $\beta $ be a
measurable mapping such that for all measurable sets $G\in\goth A$ the following estimate holds
$$
m\bigl(\beta^{-1}(G)\bigr) \le Cm(G),
$$
where the constant $C$ does not depend on $G$. For example, if the measure~$m$ is $\beta$-invariant
one can set~$C=1$. Let us consider the space $L^1(Y,m)$ of real-valued integrable functions and the
shift operator that takes every function $f\in L^1(Y,m)$ to~$f\circ\beta$. Clearly the norm of this
operator does not exceed~$C$. The mapping $\delta f :=f\circ\beta$ acts also on the space
$L^\infty(Y,m)$ and it is an endomorphism of this space. As is known, the dual space to $L^1(Y,m)$
coincides with $L^\infty(Y,m)$. Define the linear operator $A\!: L^\infty(Y,m) \to L^\infty(Y,m)$
by the identity
\begin{equation*}
\int_Y f\cdot g\circ\beta\,dm \,\equiv\, \int_Y (Af)g\,dm, \qquad g\in L^1(Y,m).
\end{equation*}

\medskip\noindent
In other words, $A$ is the adjoint operator to the shift operator in $L^1(Y,m)$. If one takes as
$g$ the index functions of measurable sets $G\subset Y$, then the latter identity takes the form
$$
\int_{\beta^{-1}(G)}f\,dm \equiv \int_G Af\,dm .
$$
Therefore $Af$ is nothing else than the Radon--Nikodim density of the additive set function
$\mu_f(G) =\int_{\beta^{-1}(G)}f\,dm$. Evidently, the operator $A$ is positive and satisfies the
homological identity
$$
A\bigl((\delta f)g\bigr) =fAg, \qquad f,g\in L^\infty(X,m).
$$

\medskip\noindent
We see that $A$ is a  transfer operator (for the $C^*$-dynamical system $(L^\infty(Y,m), \delta)$).
And in the case when $m$ is   $\beta$-invariant measure   it is a conditional expectation operator.
\end{exm}

\begin{rem}\label{2..6}
Recalling Remark \ref{2..4},\, 4)\, we have to stress that in general given a specific functional
algebra and its  endomorphism then a transfer operator is \emph{not} necessarily associated with a
\emph{point mapping} of the domain.
\end{rem}

We now present the $C^*$-dynamical analogue to Theorem \ref{1..2}.
\begin{theorem}\label{2..7}
For $\delta$-invariant functionals\/ $\mu\in M_\delta(\cal C)$ the following formula is true
\begin{equation}\label{2,,4}
 \tau(\mu)  =\inf_{n,D} \frac{1}{n} \sum_{g\in D}\mu(g) \ln\frac{\mu(A^ng)}{\mu(g)}.
 \end{equation}
\end{theorem}

\proof. This theorem can derived from  Theorem \ref{1..2}.

Indeed, to start with we observe that the Gelfand transform establishes an isomorphism between
$\cal C$ and the algebra $C(X)$ of continuous real-valued functions on the maximal ideal space $X$
of $\cal C$. Therefore  we can identify $\cal C$ with $C(X)$ mentioned above.

Moreover, under this identification the endomorphism $\delta$ mentioned in the definition of the
$C^*$-dynamical system $({\cal C},\delta)$ takes the form
\begin{equation*}
\bigl[\delta f\bigr](x) = f(\alpha (x)),
\end{equation*}
where $\alpha\!:X\to X$ is a uniquely defined continuous mapping (for details see \cite{ABL-2008},
Theorem~6.2). Thus the $C^*$-dynamical system $({\cal C},\delta)$ is completely defined by the
corresponding dynamical system $(X,\alpha)$.

In terms of the latter dynamical system the homological identity (\ref{2,,1}) for the transfer
operator $A$ can be rewritten as \eqref{1,,1}.

Since we are identifying $\cal C$ and $ C(X)$, the Riesz theorem implies that the set $M(\cal C)$
can be identified with the set $M$ of all regular Borel probability measures on~$X$ and the
identification is established by means of the formula
\begin{equation}\label{2,,5}
 \mu(\varphi) = \int_X \varphi\, d\mu,  \qquad \varphi \in{\cal C}=C(X),
\end{equation}
where $\mu$ in the right-hand part is a measure on $X$ assigned to the functional $\mu\in M (\cal
C)$ in the left-hand part.

Finally, if $\mu\in M_\delta({\cal C})$ is a $\delta$-invariant functional then the corresponding
measure $\mu$ in~\eqref{2,,5} is $\alpha$-invariant, that is
$$
\mu(f)=\mu(f\circ\alpha),\qquad f\in C(X).
$$
Therefore the set $M_\delta({\cal C})$ is naturally identified with the set $M_\alpha$ of all Borel
probability $\alpha$-in\-va\-ri\-ant measures on~$X$.

Under all these identifications the desired result follows from Theorem~\ref{1..2}. \qed


\begin{thebibliography}{9}


\bibitem{ABL-2000-1}
{A.\,B.\ Antonevich, V.\,I.\ Bakhtin, A.\,V.\ Lebedev}. Variational principle for the spectral
radius of weighted composition  and Perron--Frobenius  operators. \emph{Trudy Instituta Matematiki
NAN Belarusi} {\bf 5} (2000), 13--17 (in Russian).

\bibitem{ABL-2000-2}
{A.\,B.\ Antonevich, V.\,I.\ Bakhtin, A.\,V.\ Lebedev}. Variational principle for the spectral
radius of weighted composition  and weighted mathematical expectation operators. \emph{Doklady NAN
Belarusi} {\bf 44}(6) (2000), 7--10 (in Russian).

\bibitem{ABL-2001}
{A.\,B.\ Antonevich, V.\,I.\ Bakhtin, A.\,V.\ Lebedev}. Thermodynamics and Spectral Radius.
\emph{Nonlinear Phenomena in Complex Systems} {\bf 4}(4) (2001), 318--321.

\bibitem{ABL-2005}
{A.\,B.\ Antonevich, V.\,I.\ Bakhtin, A.\,V.\ Lebedev}. Spectra of Operators Associated with
Dynamical Systems: From Ergodicity to the Duality Principle. \emph{Twenty Years of Bialowieza: A
mathematical Anthology Aspects of Differential Geometric Methods in Physiscs}. World Scientific
Monograph Series in Mathematics. V.\ 8, Chapter 7, 129--161.

\bibitem{Bakh}
{V.\,I.\ Bakhtin}. T-entropy  and Variational principle for the spectral radius of weighted shift
operators. \emph{arXiv:0809.3106v2 [math.DS].}

\bibitem{ABL-2008} A.\,B.\ Antonevich, V.\,I.\ Bakhtin, A.\,V.\ Lebedev. T-entropy  and Variational
principle for the spectral radius of transfer and weighted shift operators. \emph{arXiv:0809.3116v2
[math.DS].}

\bibitem{Walt2}
{P.\ Walters}. \emph{An Introduction to  Ergodic Theory}. Springer-Verlag, 1982.

\end{thebibliography}
\end{document}